\documentclass{amsart}
\usepackage{amsmath}
\usepackage{amsfonts}
\usepackage{amssymb}
\usepackage{graphicx}
\theoremstyle{plain}
\newtheorem*{cpd}{Conditional Positive Definiteness}
\newtheorem*{csiqrt}{Cauchy-Schwarz Inequality for Quadratic Residue Tournaments}
\newtheorem*{cgql}{Cayley Graph Quotient Lemma}

\newtheorem*{csimn}{Cauchy-Schwarz Inequality for Manhattan Norm on $\mathbb{F}_{p^2}^n$}
\newtheorem*{sco}{Special case of $\mathbb{F}_{p}^2$}
\newtheorem*{tsi}{Triangle and Submultiplicative Inequalities in $\mathbb{F}_{p^2}$}
\theoremstyle{remark}
\newtheorem*{remark}{Remark}
\newtheorem*{remar}{Remark}
\begin{document}
\title[]{The Cauchy-Schwarz inequality in Cayley graph and tournament structures on finite fields}
\author{Stephan Foldes}
\address{Stephan Foldes \newline%
\indent Institute of Mathematics,   \newline%
\indent Tampere University of Technology,  \newline%
\indent PL 553, 33101 Tampere, Finland }
\email {sf@tut.fi}%

\author{L\'aszl\'o Major}
\address{L\'aszl\'o Major\newline%
\indent Institute of Mathematics,   \newline%
\indent Tampere University of Technology,  \newline%
\indent PL 553, 33101 Tampere, Finland }\email{laszlo.major@tut.fi}%
\hspace{-4mm} \date{Dec 24, 2009}

\begin{abstract}
The Cayley graph construction provides a natural grid structure on a finite vector space over a field of prime or prime square cardinality, where the characteristic is congruent to 3 modulo 4, in addition to the quadratic residue tournament structure on the prime subfield. Distance from the null vector in the grid graph defines a Manhattan norm. The Hermitian inner product on these spaces over finite fields behaves in some respects similarly to the real and complex case. An analogue of the Cauchy-Schwarz inequality is valid with respect to the Manhattan norm. With respect to the non-transitive order provided by the quadratic residue tournament, an analogue of the Cauchy-Schwarz inequality holds in arbitrarily large neighborhoods of the null vector, when the characteristic is an appropriate large prime.\end{abstract}
\maketitle
\section{Manhattan norms and grid graphs}
We consider the finite fields $\mathbb{F}_{p}$ and $\mathbb{F}_{p^2}$ of prime and prime square cardinality, where $p\equiv 3 \mod{4}$. The field $\mathbb{F}_{p^2}$ has a natural graph structure with the field elements as vertices, two distinct vertices $u,z$ being adjacent if $(z-u)^4=1$. The subfield $\mathbb{F}_{p}$ of $\mathbb{F}_{p^2}$ then induces a subgraph in which $x$ and $y$ are adjacent if and only if $(y-x)^2=1.$ The graph $\mathbb{F}_{p^2}$ is isomorphic to the Cartesian square $C_p^2=C_p \square C_p $, where $C_p$ is a $p$-cycle and within $\mathbb{F}_{p^2}$ the induced subgraph $\mathbb{F}_{p}$ is itself a $p$-cycle. Clearly the graph $\mathbb{F}_{p^2}$ is not planar, but can be drawn as a grid on the torus. 

For any connected graph whose vertex set is a group, the distance of any vertex $z$ from the identity element of the group is called the \emph{norm} of $z$, denoted $N(z)$. In general, distances and norms measured in connected subgraphs induced by subgroups can be larger than distances and norms measured with reference to the whole graph. However, with respect to the distance-preserving subgraph induced by $\mathbb{F}_{p}$ in $\mathbb{F}_{p^2}$, the norm of any $z\in \mathbb{F}_{p}$ is the same as its norm with respect to the whole graph $\mathbb{F}_{p^2}$: this is simply the length of the shortest path from $0$ to $z$ in the cycle induced by $\mathbb{F}_{p}$.

For $q=p$ or $q=p^2$, the $n$-dimensional vector space $\mathbb{F}_{q}^n$ is also endowed with the Cartesian product graph structure $\mathbb{F}_{q}\square \cdots \square \mathbb{F}_{q}$ isomorphic to $C_p^n$ or $C_p^{2n}$. The norm of a vector $\textbf{v}=(v_1,\ldots,v_n)$ in $\mathbb{F}_{q}^n$ is then equal to the sum $N(v_1)+\cdots + N(v_n)$ and we also write $N(\textbf{v})$ for this vector norm. 

The Gaussian integers $\mathbb{Z}[i]$ also constitute a graph in which $u$ and $z$ are adjacent if and only if $(z-u)^4=1$. 

It is easy to see that the norm in this \emph{infinite Manhattan grid} satisfies the triangle and submultiplicative inequalities
\begin{align*}
	N(u+z)&\leq N(u)+ N(z)\\
	N(uz)&\leq N(u)N(z)
\end{align*}

To emphasize that the norms on $\mathbb{F}_{p^2}$, $\mathbb{F}_{p^2}^n$ and $\mathbb{Z}[i]$   
are understood with reference to the specific grid graphs defined above, we call these norms \emph{Manhattan norms}.
Throughout this paper we thi\vspace{2mm}nk of $\mathbb{F}_{p^2}$ as the ring quotient $\mathbb{Z}[i]/(p).$

\section{Graph quotients and Cayley graphs}
Given a graph $G$ (undirected, with possible loops) on vertex set $V$ and an equi-valence relation $\equiv$ on $V$, the \emph{quotient graph} $G/\hspace{-0.7mm}\equiv$ is defined as follows: the vertices of $G/\hspace{-0.85mm}\equiv$ are the equivalence classes of $\equiv$, and classes $A, B$ are adjacent if for some $a\in A$, $b\in B$, the elements $a,b$ are adjacent in $G$. Note that the distance of $A$ to $B$ in the quotient graph is at most equal to, but possibly less than the minimum of the distances $a$ to $b$ for all $a\in A$, $b\in B$. Note also that $G/\hspace{-0.75mm}\equiv$ can have loops even if $G$ has not.

Given a group $G$ with identity element $e$ and a set $\Gamma$ of group elements that generates $G$, the \emph{(left) Cayley graph} $\mathcal{C}(G,\Gamma)$ of $G$ with respect to $\Gamma$ has vertex set $G$, elements $a,b\in G$ being considered adjacent if $ab^{-1}$ or $ba^{-1}$ belongs to $\Gamma$. For each congruence $\equiv$ of the group $G$, corresponding to some normal subgroup $H$, $\Gamma$ yields a generating set $\Gamma_{\equiv}$ of $G/\hspace{-0.7mm}\equiv$ consisting with those classes of $\equiv$ that intersect $\Gamma$. The graph quotient of $\mathcal{C}(G,\Gamma)$ by the equivalence $\equiv$ coincides with the Cayley graph of the quotient graph $G/\hspace{-0.7mm}\equiv$ with respect to $\Gamma_{\equiv}$. For $R\subseteq G$ inducing a connected subgraph $[R]$ in $\mathcal{C}(G,\Gamma)$, denote by $d_R(x,y)$ the distance function of the subgraph $[R]$. Denoting by $xH$ the $H$-coset of any $x \in G$, this relates to norms in $\mathcal{C}(G,\Gamma)$ and $\mathcal{C}(G,\Gamma)/\hspace{-0.7mm}\equiv$ as follows: for all $x\in R$,
$$ d_R(x,e)\vspace{-0mm} \geq N(x)\geq N(xH)$$
Both inequalities can be strict. \vspace{1.5mm}
However, we have:
\begin{cgql}Let a group $G$ with identity $e$ be generated by $\Gamma \subseteq G$, and consider any normal subgroup $H$ with corresponding congruence $\equiv$. There is a set $R\subseteq G$ having exactly one element in common with each congruence class modulo $H$, and such that for every $x \in R$ $$d_R(x,e) =  N(x) =  N(xH)$$\end{cgql}
\begin{proof}We can define the unique (representative) element $r(A)\in R\cap A$ for each coset $A$ by induction on the distance $d(H,A)$ of $A$ from $H$ in $\mathcal{C}(G,\Gamma)/\hspace{-0.7mm}\equiv$. Let $r(H)=e$. Assuming $r(A)$ defined for all $A$ with $d(H,A)\leq m$, let a coset $B$ have distance $m+1$ from $H$. Choose any coset $A$ adjacent to $B$ with $d(H,A)=m$ and elements $a\in A$, $b\in B$ that are adjacent in $\mathcal{C}(G,\Gamma)$. Let $r(B)=ba^{-1}r(A)$.\end{proof}

We can apply the above lemma in the case where  $G=\mathbb{Z}[i]$, $\Gamma=\{1,i\}$ and $H=p\mathbb{Z}[i]=\{pa+pbi: a,b\in \mathbb{Z}\}$ for a prime integer $p\equiv 3 \mod{4}$. Now $\mathcal{C}(G,\Gamma)$ and $\mathcal{C}(G,\Gamma)/\hspace{-0.7mm}\equiv$ are the Manhattan grid graphs on $\mathbb{Z}[i]$ and $\mathbb{Z}[i]/H=\mathbb{F}_{p^2}$, respectively. Referring to the set $R$ of representatives in the lemma, for any $H$-cosets $X,Y$ let $x,y$ be the unique elements in $X\cap R$, $Y\cap R$. As $xy\in XY$, we have $N(XY)\leq N(xy)$. By the submultiplicative inequality in $\mathbb{Z}[i]$ we have $N(xy)\leq N(x)N(y)$. Using the lemma we have $N(x)N(y)=N(X)N(Y)$. This yields a submultiplicative inequality in $\mathbb{F}_{p^2}$ and a similar reasoning on the coset $X+Y$ yields a triangle inequality: 

\begin{tsi} For all $u, z$ in $\mathbb{F}_{p^2}$
\begin{align*}
	N(u+z)&\leq N(u)+ N(z)\\
	N(uz)&\leq N(u)N(z) 
\end{align*}\hspace{124.14mm}$\square$
\end{tsi}

This indicates that Manhattan distance provides a well-behaved notion of neighborhood of $0$ in the finite fields $\mathbb{F}_{p^2}$.

\section{Squares in $\mathbb{F}_p$ and non-transitive order}
For each prime $p\equiv 3 \mod{4}$ the \emph{quadratic residue tournament} on $\mathbb{F}_{p}$ is the directed graph with vertex set $\mathbb{F}_{p}$ in which there is an arrow from vertex $x$ to vertex $y$ if $y-x$ is a non-zero square in $\mathbb{F}_{p}$, in which case we write $x<_p y$. We write $x \leq_p y$ if $x<_p y$ or $x=y$. The relation $\leq_p$ is reflexive, anti-symmetric but not transitive, and for every $x\neq y$ exactly one of $x \leq_p y$ or $y \leq_p x$ holds. Using Dirichlet's theorem on primes in arithmetic progressions, Kustaanheimo showed \cite{K1} that for every positive integer $k$, there is a prime $p\equiv 3 \mod{4}$, such that $\leq_p$ is a transitive (and linear) order relation on $\{0,1,\ldots,k\}\subseteq \mathbb{F}_p$, that is, all positive integers up to $k$ are quadratic residues$\mod{p}$. Obviously $k$ cannot exceed $(p-1)/2$. Implications of \cite{K1} and related questions were investigated by  J\"{a}rnefelt, Kustaanheimo, Quist \cite{J1,Q1}, in particular with a view to discrete models in physics, also in subsequent application-oriented work between the 1950's (Coish \cite{C1}) and the 1980's (Nambu \cite{N1}). For further references see \cite{F1}. In particular \cite{K1} implies that for every positive integer $k$, there is a prime $p\equiv 3 \mod{4}$, such that all $z\in \mathbb{F}_{p^2}$ with $N(z)\leq k$ are squares in $\mathbb{F}_{p^2}$. (Note that all elements of the prime subfield $\mathbb{F}_{p}$ are squares in $\mathbb{F}_{p^2}$.) To emphasise the analogy of the relation $\leq_p$ with the ordinary inequality relation $\leq$ among numbers, we say that a non-zero $z\in \mathbb{F}_{p^2}$ is \emph{positive} if $z\in \mathbb{F}_{p} $ and $0\leq_p z$.

\section{Inner products compared in non-transitive order}
The only non-trivial automorphism of the field $\mathbb{F}_{p^2}$ associates to each $z\in \mathbb{F}_{p^2}$ its \emph{conjugate} $\overline{z}$. The \emph{inner product} $\textbf{v}\cdot \textbf{w}$ of vectors $\textbf{v}=(v_1,\ldots,v_n)$ and $\textbf{w}=(w_1,\ldots,w_n)$ in $\mathbb{F}_{p^2}^n$ is defined as the scalar $v_1\overline{w_1}+\cdots + v_n \overline{w_n} \in \mathbb{F}_{p^2}$. This inner product is left and right distributive over vector addition, satisfies $\textbf{v}\cdot \textbf{w}=\overline{\textbf{w}\cdot \textbf{v}}$, $c(\textbf{v}\cdot \textbf{w})=(c\textbf{v})\cdot \textbf{w}=\textbf{v}\cdot (\overline{c}\textbf{w})$ for all $c\in \mathbb{F}_{p^2}$. However, while $\textbf{v}\cdot \textbf{v}$ belongs to the prime subfield $\mathbb{F}_{p}$, $\textbf{v}\cdot \textbf{v}$ is not necessarily positive, and can be $0$ even if $\textbf{v}\neq \textbf{0}$. Still, a conditional version of positive definiteness holds locally:

\begin{cpd}
\upshape 
 \textit{For every $k\geq 1$ there is a prime $p\equiv 3 \mod{4}$, such that for all $n\geq 1$ and for all vectors} $\textbf{v}\in \mathbb{F}_{p^2}^n$ \textit{of Manhattan norm }$N(\textbf{v})\leq k$, \textit{we have} $0\leq_p \textbf{v}\cdot \textbf{v}$ \textit{with equality if and only if} $\textbf{v}= \textbf{0}$.
\end{cpd}
\begin{proof}By Kustaanheimo's result in \cite{K1} there is a prime integer $p\equiv 3 \mod{4}$ such that $0,1,\ldots,2k^3$ are all quadratic residues mod $p$. For  $\textbf{v}=(v_1,\ldots,v_n)$ in $\mathbb{F}_{p^2}^n$, let $v_j=a_j+b_ji$, where $i^2=-1$. If $N(\textbf{v})\leq k$ then for all $j$, $N(a_j)\leq k$ and $N(b_j)\leq k$, $v_j\overline{v_j}=a_j^2+b_j^2$ belongs to the set of squares $\{0,\ldots,2k^2\}$. Since $v_j$ can be non-zero for at most $k$ indices $1\leq j\leq n$ only, the sum of the corresponding terms $a_j^2+b_j^2$
belongs to the set of squares $\{0,1,\ldots,2k^3\}.$ \end{proof}
 Note that for all vectors $\textbf{v},\textbf{w}\in \mathbb{F}_{p^2}^n$
 \begin{align*}
	(\textbf{v}\cdot \textbf{w})(\textbf{w}\cdot \textbf{v})=(\textbf{v}\cdot \textbf{w})(\overline{\textbf{v}\cdot \textbf{w}})&\in \mathbb{F}_{p}\hspace{2mm}\text{and}\\  (\textbf{v}\cdot \textbf{v})(\textbf{w}\cdot \textbf{w})&\in \mathbb{F}_{p} \end{align*}
If $\textbf{v}$ and $\textbf{w}$ are \emph{proportional}, i.e. if there exists a scalar $c$ in $\mathbb{F}_{p^2}$ such that $\textbf{v}=c\textbf{w}$ or $\textbf{w}=c\textbf{v}$, then the above two products are equal. Generally, they are related in the quadratic residue tournament of $\mathbb{F}_{p}$ as follows.

\begin{csiqrt}\upshape 
\textit{For eve-ry $k\geq 1$ there is a prime $p\equiv 3 \mod{4}$, such that for all $n\geq 1$ and for all vectors }$\textbf{v},\textbf{w}\in \mathbb{F}_{p^2}^n$ \textit{of Manhattan norm at most} $k$, $$(\textbf{v}\cdot \textbf{w})(\textbf{w}\cdot \textbf{v}) \leq_p (\textbf{v}\cdot \textbf{v})(\textbf{w}\cdot \textbf{w})$$\end{csiqrt}
\begin{proof}For $n=1$ the inequality holds trivially as the two sides are equal. Assume $n\geq 2$, $\textbf{v}=(v_1,\ldots,v_n),\textbf{w}=(w_1,\ldots,w_n)$. For all $1\leq i\leq n$, $N(v_i)\leq k$, $N(w_i)\leq k$. By Kustaanheimo's result \cite{K1} there is a prime $p\equiv 3 \mod{4}$ such that all positive integers up to $4k^6$ are quadratic residues modulo $p$. For each of the $\binom n2$ pairs $\{i,j\}\subseteq \{1,\ldots,n\}$, $i\neq j$, by the triangle and submultiplicative inequalities in $\mathbb{F}_{p^2}$
$$N[(v_iw_j-v_jw_i)\vspace{-2mm}(\overline{v}_i\overline{w}_j-\overline{v}_j\overline{w}_j)]\leq (k^2+k^2)^2=4k^4$$
Thus the element $$(v_iw_j\overline{v}_i\overline{w}_j+v_jw_i\overline{v}_j\overline{w}_i)-(v_iw_j\overline{v}_j\overline{w}_i+v_jw_i\overline{v}_i\overline{w}_j)=(v_iw_j-v_jw_i)(\overline{v}_i\overline{w}_j-\overline{v}_j\overline{w}_j)$$ is a square of Manhattan norm at most $4k^4$ in $\mathbb{F}_{p}$, and it is non-zero for at most $\binom k2 \leq k^2$ pairs $\{i,j\}$.
Summing over all pairs $\{i,j\}$, all but at most $\binom k2 \leq k^2$ terms vanish in the sum $$\sum[(v_iw_j\overline{v}_i\overline{w}_j+v_jw_i\overline{v}_j\overline{w}_i)-(v_iw_j\overline{v}_j\overline{w}_i+v_jw_i\overline{v}_i\overline{w}_j)]$$ which therefore has Manhattan norm at most $4k^6$ and it must also be a square in $\mathbb{F}_{p}$. But this sum is equal to the difference of products
$$\sum_{i=1}^nv_i\overline{v}_i\sum_{j=1}^nw_j\overline{w}_j-\sum_{i=1}^nv_i\overline{w}_i\sum_{j=1}^n\overline{v}_jw_j=(\textbf{v}\cdot \textbf{v})(\textbf{w}\cdot \textbf{w})-(\textbf{v}\cdot \textbf{w})(\textbf{w} \cdot \textbf{v})$$ which is consequently a square in $\mathbb{F}_{p}$.\end{proof}
\begin{remark}From the proof it is clear that, in analogy with the classical Cauchy-Schwarz inequality, for vectors $\textbf{v},\textbf{w}$ of norm not exceeding $k$ in $\mathbb{F}_{p^2}^n$, where $p$ is related to $k$ as stipulated above, the Cauchy-Schwarz inequality with respect to $\leq_p$ holds with equality if and only if $v_iw_j-v_jw_i=0$ for all $i,j$, i.e. if and only if $\textbf{v},\textbf{w}$ are proportional. \end{remark}

We note that the inequality established above is conditional, it holds only in a specified Manhattan neighborhood of the null vector. Every non-zero element of $\mathbb{F}_{p}$ can be written as a sum of two squares, in particular there are $a,b\in \mathbb{F}_{p}$, such that $a^2+b^2=-1$. For $z=a+bi$ we have $z \overline{z}=-1$. As soon as $n\geq 2$, in $\mathbb{F}_{p^2}^n$ let $$\textbf{v}=(a,b,0,\ldots,0)\hspace{3mm}\text{and}\hspace{3mm}\textbf{w}=(bz,-az,0,\ldots,0)$$
The inequality $(\textbf{v}\cdot \textbf{w})(\textbf{w}\cdot \textbf{v}) \leq_p (\textbf{v}\cdot \textbf{v})(\textbf{w}\cdot \textbf{w})$ fails because the left-hand side is $0$ and the right-hand side is $-1$. In fact if $n\geq 3$, the inequality can be invalidated with vectors $\textbf{v},\textbf{w}$ in $\mathbb{F}_{p}^n$ as follows. Taking again $a,b\in  \mathbb{F}_{p}$ with $a^2+b^2=-1$, let
$$\textbf{v}=(1,a,b,0,\ldots,0)\hspace{3mm}\text{and}\hspace{3mm}\textbf{w}=(1,0,0,0,\ldots,0)$$
However, the Cauchy-Schwarz inequality holds unconditionally in the $2$-dimensional case for vectors with components in $\mathbb{F}_{p}$:

\begin{sco}\upshape \textit{Let $p$ be a prime congruent $3$ modulo $4$. For all vectors} $\textbf{v},\textbf{w}$ \textit{in $\mathbb{F}_{p}^2$}
$$(\textbf{v}\cdot \textbf{w})(\textbf{w}\cdot \textbf{v}) \leq_p (\textbf{v}\cdot \textbf{v})(\textbf{w}\cdot \textbf{w})$$\end{sco}

\begin{proof}Now the conjugation appearing in the inner products is the identity. Written in components,
\begin{align*}
(\textbf{v}\cdot \textbf{v})(\textbf{w}\cdot \textbf{w})-(\textbf{v}\cdot \textbf{w})(\textbf{w}\cdot \textbf{v})&=(v_1^2+v_2^2)(w_1^2+w_2^2)-(v_1w_1+v_2w_2)^2=\\
=v_1^2w_2^2+v_2^2w_1^2-2v_1w_1v_2w_2&=(v_1w_2-v_2w_1)^2	
\end{align*}
\end{proof}
\section{Manhattan norm of inner product}
The Manhattan norm can be seen to be submultiplicative not only on the ring $\mathbb{Z}[i]$ and its quotient field $\mathbb{F}_{p^2}$, but on all vector spaces $\mathbb{F}_{p^2}^n$, with respect to the inner product:
\begin{csimn}\upshape \textit{Consider any prime $p \equiv 3 \mod{4}$ and let $n\geq 1$. For all} $\textbf{v},\textbf{w}\in \mathbb{F}_{p^2}^n$ $$N(\textbf{v}\cdot \textbf{w})\leq N(\textbf{v}) N(\textbf{w})$$
\end{csimn}
\begin{proof}Let $\textbf{v}=(v_1,\ldots,v_n),\textbf{w}=(w_1,\ldots,w_n)\in \mathbb{F}_{p^2}^n$. Then $\textbf{v}\cdot \textbf{w}=\sum v_j\overline{w_j}$. Clearly $N(z)=N(\overline{z})$ for any $z\in \mathbb{F}_{p^2}$. By the triangle  and submultiplicative inequalities in  $\mathbb{F}_{p^2}$ we have 
\begin{align*}
N(\textbf{v}\cdot \textbf{w})= N\big(\textstyle\sum{ v_j\overline{w_j}}\big)&\leq \textstyle\sum N\big({ v_j\overline{w_j}}\big)\leq \sum N(v_j)N(w_j)\leq\\&\leq \textstyle\sum N(v_j) \sum N(w_j)= 	N(\textbf{v}) N(\textbf{w})
\end{align*}\end{proof}
\begin{remar}The inequality $N(\textbf{v}\cdot \textbf{w})\leq N(\textbf{v}) N(\textbf{w})$ is easily interpreted and continues to hold for $\textbf{v},\textbf{w}$ in the module $(\mathbb{Z}[i]/m\mathbb{Z}[i])^n$ for any positive integer $m$. As soon as $m$ is composite, or a prime not congruent to 3 modulo 4, the ring $\mathbb{Z}[i]/m\mathbb{Z}[i]$ fails to be an integral domain.
\end{remar}

\end{document}